\documentclass{siamart0516}
\usepackage{lipsum}
\usepackage{amssymb}
\usepackage{amsmath}
\usepackage{amsfonts}
\usepackage{graphicx}
\usepackage{epstopdf}
\usepackage{algorithmic}
\ifpdf
  \DeclareGraphicsExtensions{.eps,.pdf,.png,.jpg}
\else
  \DeclareGraphicsExtensions{.eps}
\fi
\numberwithin{equation}{section}
\usepackage{tcolorbox}

\topmargin by-\headsep  \textheight   236mm   
\def\beqn{\begin{eqnarray}}  \def\eqn{\end{eqnarray}}
\def\beqnx{\begin{eqnarray*}} \def\eqnx{\end{eqnarray*}}
\def\bn{\begin{enumerate}} \def\en{\end{enumerate}}

\newcommand{\TheTitle}{Adaptive Residual-Driven Newton Solver} 
\newcommand{\TheAuthors}{Renjie Ding and Dongling Wang}
\headers{\TheTitle}{\TheAuthors}

\title{Adaptive Residual-Driven Newton Solver for Nonlinear Systems of Equations}

\author{
Renjie Ding\thanks{Hunan Key Laboratory for Computation and Simulation in Science and Engineering, School of Mathematics and Computational Science,  Xiangtan University, Xiangtan, Hunan 411105, P.R. China. (\email{drjmath@smail.xtu.edu.cn}).}
\and
 Dongling Wang\thanks{Hunan Key Laboratory for Computation and Simulation in Science and Engineering, School of Mathematics and Computational Science,  Xiangtan University, Xiangtan, Hunan 411105, P.R. China.
The work of D. Wang was partially supported by National Natural Science Foundation of China under grants 12271463. 
(\email{wdymath@xtu.edu.cn}).}
  }
 
\ifpdf
\hypersetup{
  pdftitle={\TheTitle},
  pdfauthor={\TheAuthors}
}
\fi

\begin{document}

\maketitle

\begin{abstract}
Newton-type solvers have been extensively employed for solving a variety of nonlinear system of algebraic equations. However, for some complex nonlinear system of algebraic equations, efficiently solving these systems remains a challenging task. The primary reason for this challenge arises from the unbalanced nonlinearities within the nonlinear system. Therefore, accurately identifying and balancing the unbalanced nonlinearities in the system is essential. In this work, we propose a residual-driven adaptive strategy to identify and balance the nonlinearities in the system. The fundamental idea behind this strategy is to assign an adaptive weight multiplier to each component of the nonlinear system, with these weight multipliers increasing according to a specific update rule as the residual components increase, thereby enabling the Newton-type solver to select a more appropriate step length, ensuring that each component in the nonlinear system experiences sufficient reduction rather than competing against each other. More importantly, our strategy yields negligible additional computational overhead and can be seamlessly integrated with other Newton-type solvers, contributing to the improvement of their efficiency and robustness. 
We test our algorithm on a variety of benchmark problems, including a chemical equilibrium system, a convective diffusion problem, and a series of challenging nonlinear systems. The experimental results demonstrate that our algorithm not only outperforms existing Newton-type solvers in terms of computational efficiency but also exhibits superior robustness, particularly in handling systems with highly imbalanced nonlinearities.
\end{abstract}



\section{Introduction}

The numerical solution of nonlinear system of algebraic equations has been one of the most extensively studied topics in computational mathematics. Many applications, including the discretization of nonlinear partial differential equations (PDEs), the implicit Runge-Kutta method, and numerous engineering problems, require the solution of nonlinear systems of algebraic equations. Classical Newton-type methods, such as the Jacobian-Free Newton-Krylov (JFNK) technique \cite{knoll2004jacobian} and inexact Newton methods \cite{dennis1996numerical}, are commonly employed for solving these nonlinear systems. While these methods exhibit remarkable effectiveness in many cases, they can sometimes encounter challenges. Specifically, they may suffer from slow convergence or even fail to converge altogether, particularly when a suitable initial guess is lacking or when the nonlinear components within the system are unbalanced.

To accelerate the convergence of Newton-type methods, researchers have devised various strategies, which can broadly be categorized into three groups: improvements in line search \cite{knoll2004jacobian, brown1990hybrid, pawlowski2006globalization}, employing nonlinear preconditioning techniques \cite{cai2002nonlinearly, cai2011inexact, klawonn2014nonlinear, tang2019fully, luo2023preconditioned, huang2016nonlinearly, hwang2015parallel}, and establishing a suitable initial guess \cite{kim2006newton, huang2020int, choi2022choice, kim2017newton}. The line search strategy is beneficial for globalizing inexact Newton methods, but the step length obtained from line searches are often determined by strongly nonlinear components. This can lead to stagnation in the nonlinear residual curve \cite{gropp2000globalized, cai2002nonlinearly}. Nonlinear preconditioning is a technique aimed at improving the robustness and efficiency of Newton-type solvers by balancing the system by eliminating local high nonlinearities that can hinder the performance of the solver in nonlinear systems. Despite the fact that many preconditioning techniques are already quite effective, there are still challenges that need to be addressed. For example, in the region-based approach \cite{luo2019nonlinear}, the number of slow components that need elimination heavily relies on the selection of preselected parameters, and an inappropriate choice can significantly compromise the effectiveness of the preconditioner. The advantages of selecting a good initial guess for any iterative algorithm are evident, yet it often poses a significant challenge. Generally, finding a universal algorithm to choose a good initial guess is difficult and typically requires additional analysis tailored to the specific problem. For time-varying problems, one possible approach is to use the result from the previous time step or leverage the linearization of the equation at the current time step \cite{choi2022choice}. For other types of equations, considering the solution of a simpler equation as the initial value for a more complex one can be helpful \cite{kim2017newton}. For instance, in the case of compressible flow problems at high Reynolds numbers, the solution of a problem at low Reynolds numbers can be selected as the initial guess.

In this work, we propose a residual-driven adaptive weighting strategy that dynamically adjusts the weights of different components in the nonlinear system 
during the iteration process to address unbalanced nonlinear equations, thereby accelerating the convergence of Newton-type methods. Our improvement strategy is both concise and efficient, significantly enhancing the robustness of the algorithm with virtually no additional computational cost. Furthermore, the algorithm is highly versatile and can be seamlessly integrated with other acceleration strategies to further expedite the convergence of Newton-type methods.

The paper is structured as follows. In Section \ref{sec:New}, we provide a review of the inexact Newton method with backtracking and briefly analyze the factors contributing to its slow convergence. In Section \ref{sec:ARDN}, we introduce the proposed adaptive residual-driven Newton-type solver.
In Section \ref{sec:addopt}, we also discuss other techniques designed to enhance the convergence of Newton-type solvers, comparing them with our method.
Furthermore, it is shown that our method can be effectively integrated with other acceleration strategies.
In Section \ref{sec:exam}, we present numerical experiments conducted on our algorithm, including tests on the chemical equilibrium system, the convective diffusion problem, and a variety of challenging problems. These experiments are intended to validate the effectiveness of our proposed algorithm, assess its robustness and efficiency, and compare its performance with other Newton-type solvers. Finally, Section \ref{sec:con} offers some concluding remarks and key observations.

\section{Newton's iteration and its improvements}
\label{sec:New}
Consider the following nonlinear system of algebraic equations
\begin{equation}\label{equation:2.1}
F(X)=0,
\end{equation}
where $F:\mathbb{R}^n \rightarrow \mathbb{R}^n, n\in \mathbb{N}^+$ is continuously differentiable, and we denote $F=(F_1,\dots,F_n)^T$, where $F_i=F_i(X)$ and $X=(X_1, \dots, X_n)^T$. Now, we aim to seek a vector $X^* \in \mathbb{R}^n$ such that $F(X^*)=0$, starting from an initial guess $X^0 \in \mathbb{R}^n$. The classic method for solving nonlinear system of algebraic equations \eqref{equation:2.1} is Newton's method \cite{yamamoto2000historical}.  Its main idea is local linearization. Specifically, if we have the current approximate solution $X^k$, the next approximate solution $X^{k+1}$ is determined by the following formula:
\begin{equation}\label{equation:2.2}
X^{k+1}=X^k+S^k,
\end{equation}
where $S^k$ is referred to as the Newton direction, which is determined by the following Newton's equation:
\begin{equation}\label{equation:2.3}
F^{\prime}(X^k)S^k=-F(X^k).
\end{equation}
The key advantage of the Newton method is that it exhibits local quadratic convergence  if the Jacobian $F^{\prime}$ is Lipschitz continuous near $X^*$ \cite{ortega2000iterative}. However, its drawback lies in the necessity to solve a linear system \eqref{equation:2.3} at each Newton iteration. When $n$ is large, the cost of solving this linear system \eqref{equation:2.3} becomes excessively expensive. Additionally, Newton's method lacks global convergence properties.

In order to overcome the high computational cost associated with solving linear systems, the inexact Newton condition \cite{dembo1982inexact} was introduced as an alternative to Newton's equation \eqref{equation:2.3} for obtaining $S^k$, which is then referred to as the inexact Newton direction. Specifically, the inexact Newton direction $S^k$ needs to satisfy the following inexact Newton condition:
\begin{equation}\label{equation:2.4}
\Vert F^{\prime}(X^k)S^k+F(X^k)\Vert\leq \eta^{k}\Vert F(X^k)\Vert,
\end{equation}
where $\eta^k\in(0,1)$ is referred to as the forcing term, which governs the accuracy of solving the linear system \eqref{equation:2.3}. The choice of $\eta^k$ reflects the underlying idea of the inexact Newton condition. Generally, when the current approximate solution $X^k$ is far from the exact solution $X^*$, solving the linear system accurately is unnecessary, and a larger value of $\eta^k$ should be chosen. Conversely, a smaller value should be taken when the approximate solution $X^k$ is close to the exact solution $X^*$. Currently, the selection methods proposed by Eisenstat and Walker \cite{eisenstat1996choosing} are the most influential, which is given by
\begin{equation}\label{equation:2.5}
\eta^{k}=\left\{\begin{array}{ll}\eta_0,&\|F(X^{k})\|\geq\beta,\\[1ex]\frac{\left|\|F(X^{k})\|-\|F^{'}(X^{k-1})S^{k-1}+F(X^{k-1}))\|\right|}{\|F(X^{k-1})\|},&\|F(X^{k})\|<\beta,\end{array}\right.
\end{equation}
where $\eta_0$ and $\beta$ are two given constants. In addition, to solve the linear system \eqref{equation:2.3} based on the inexact condition \eqref{equation:2.4}, we typically use Krylov subspace iterative methods \cite{saad2003iterative}, and we refer to solving the linear system \eqref{equation:2.3} based on \eqref{equation:2.4} as solving the linear system at level $\eta^k$. Furthermore, we should note that for some strongly nonlinear system of algebraic equation, the resulting linear systems are often difficult to solve. In such cases, we should construct appropriate preconditions based on the specific linear system.

In order to improve the local convergence of Newton's method, it is necessary to adopt globalization techniques. Many scholars have proposed strategies to enhance Newton's method, leading to several globally convergent methods \cite{kelley1995iterative, eisenstat1994globally, bellavia2001globally, dennis1996numerical, an2007globally}. 
Here we recall
the basic idea of line search method with a backtracking strategy \cite{dennis1996numerical}. 
For the equation \eqref{equation:2.1},  the function for performing the line search is as follows:
\begin{equation}\label{equation:2.6}
f(X)=\frac{\Vert F(X)\Vert^2}{2}=\frac{1}{2}\sum_{i=1}^{n}F_i(X)^2,
\end{equation}
where $f$ is referred to as the merit function. 
We now need to determine the step length $\lambda^k$ along the direction $S^k$ using a line search criterion. 
And we can use the following Armijo conditions \cite{wolfe1969convergence} to determine $\lambda^k$ along the inexact Newton direction $S^k$ as
\begin{equation}\label{equation:2.7}
f(X^k+\lambda^kS^k)\leq f(X^k)+\alpha \lambda^{k}\nabla f(X^k)^TS^k,
\end{equation} 
where the parameter $\alpha$ is used to ensure sufficient decrease in the merit function $f$, and it is typically set to $\alpha=10^{-4}$. 
Furthermore, we typically set the maximum number of line search iterations $g_{\max}$ and the line search step length decay rate $\rho$. 
Therefore, we have obtained an inexact Newton method with backtracking \cite{dennis1996numerical}. 
Assuming $X^k$ is the current approximate solution, the new approximate solution $X^{k+1}$ is given by
\begin{equation}\label{equation:2.8}
X^{k+1}=X^k+\lambda^{k}S^k,
\end{equation} 
where the inexact Newton direction $S^k$ and the step length $\lambda^k$ are determined by equations \eqref{equation:2.4} and \eqref{equation:2.7}, respectively. Additionally, the stopping criterion for the nonlinear iteration is given by
\begin{equation}\label{equation:2.9}
\Vert F(X^k)\Vert\leq \max \left\{ \gamma_a, \gamma_r \Vert F(X^0)\Vert \right\},
\end{equation} 
where $\gamma_a$ and $\gamma_r$ are typically referred to as the absolute tolerance and the relative tolerance, respectively.

The Inexact Newton method with backtracking (INB) has several advantages over the classical Newton method, primarily in the following two aspects: (a) The INB obtains the inexact Newton direction $S^k$ by satisfying the inexact Newton condition (\ref{equation:2.4}), which typically results in a lower computational cost. (b) By introducing the backtracking line search criterion, the INB ensures global convergence \cite{kelley1995iterative, eisenstat1994globally, bellavia2001globally, dennis1996numerical, an2007globally}. 

Although the INB is quite effective, its convergence can still be very slow for nonlinear systems with unbalanced nonlinearities. Specifically, it is important to note that the step length $\lambda^k$ is a crucial parameter in the INB. When the value of $\lambda^k$ is very small, it often leads to slow convergence of the INB. The selection of the step length is related to the Armijo conditions (\ref{equation:2.7}) used in the line search. The range of step length determined by the Armijo conditions is influenced by the term $\alpha \nabla f(X^k)^TS^k$. If the linear system (\ref{equation:2.3}) is solved exactly, then the expression for $\alpha \nabla f(X^k)^TS^k$ is as follows:
\begin{equation}\label{equation:2.10}
\alpha \nabla f(X^k)^TS^k=-\alpha \nabla f(X^k)^TF^{\prime}(X^k)^{-1}F(X^k)=-\alpha \Vert F(X^k)\Vert^2=-\alpha\sum_{i=1}^nF_i(X^k)^2.
\end{equation} 
Usually, a larger value of $|\alpha \nabla f(X^k)^T S^k|$ results in a smaller step length selected through the line search, and vice versa. From equation (\ref{equation:2.10}), it is known that the value of the step length $\lambda^k$ is determined by the components that contribute significantly to the residual norm in the nonlinear system.

In order to compare with Algorithm \ref{alg:ARDN} that we will be developing below, we summarize the above ideas for INB as Algorithm \ref{alg:INM}.

\begin{algorithm}
\caption{INB: The Inexact Newton Method with Backtracking}
\label{alg:INM}
\begin{algorithmic}[1]
\STATE \textbf{Input:} $X^0$, $\eta_0$, $\beta$, $\gamma_a$, $\gamma_r$, $\alpha$, $\rho$
\STATE \textbf{Output:} $X^*$
\STATE $k\leftarrow 0$ \hfill  $\triangleright$ Nonlinear iteration index
\STATE $f(X):=\frac{\Vert F(X)\Vert^2}{2}$ \hfill$\triangleright$ Define the merit function
\WHILE{$\Vert F(X^k)\Vert> \max\{\gamma_r\Vert F(X^0)\Vert, \gamma_a\}$}
\STATE // Compute the forcing term $\eta^k$ using equation (\ref{equation:2.5}).
\STATE Select $\eta^k\in[0, 1)$\
\STATE // Precondtioner $\mathcal{P}^k$
\STATE $S^k \leftarrow \mathrm{Krylov}(F^{\prime}(X^k),-F(X^k), \mathcal{P}^k)$\hfill$\triangleright$ Solve the linear system (\ref{equation:2.3}) (at the level $\eta^k$)
\FOR{$i \in \{0, \dots, g_{max}\}$}
    \STATE // Obtain the step length through backtracking line search.
    \IF{$f(X^k + \lambda_i S^k) \leq f(X^k) + \alpha \lambda_i \nabla f(X^k)^T S^k$}
        \STATE $\lambda^k\leftarrow \lambda_i$
        \STATE break
        \ELSE
        \STATE $\lambda_i \leftarrow \rho \lambda_i $
    \ENDIF
\ENDFOR
\STATE $k\leftarrow k+1$ \hfill  $\triangleright$ Update nonlinear iteration index
\STATE $X^k \leftarrow X^k+\lambda^kS^k$   \hfill  $\triangleright$ Update the approximate solution
\ENDWHILE
\STATE $X^*\leftarrow X^k$ 
\end{algorithmic}
\end{algorithm}

\section{Adaptive Residual-Driven Newton Solver}
\label{sec:ARDN}
To improve the inexact Newton method with backtracking, we now introduce a new adaptive residual-driven Newton solver, which is refer to ARDN.

\subsection{ARDN: Adaptive Residual-Driven Newton Solver}
The goal of designing ARDN is to accurately identify and balance the unbalanced nonlinearities in the nonlinear system of algebraic equations. More specifically, we aim to design a simple and low-cost scheme that enables the Newton-type solver to select a more appropriate step length at each iteration, and to prioritize the elimination of components that contribute significantly to the residual norm in the current iteration, combining this with equation (\ref{equation:2.10}), this approach not only ensures a sufficient decrease in the merit function (\ref{equation:2.6}), but also balances the unbalanced nonlinearities in the system, thereby facilitating subsequent Newton iterations.

To achieve this goal, for each component $F_i$ for $i=1,\dots, n$ in the nonlinear  system of algebraic equations, we assign an adaptive weight multiplier $\omega_i^k$, where the index $k$ denotes the $k$-th iteration. These weight multipliers are updated along with the Newton iterations progress. Therefore, the form of the merit function (\ref{equation:2.6}) will be different at each Newton iteration. The merit function $f^k(X)$ at the $k$-th iteration is given by 
\begin{equation}\label{equation:3.1}
f^k(X)=\frac{\Vert \omega^k\odot F(X)\Vert^2}{2}=\frac{1}{2}\sum_{i=1}^{n}(\omega_i^kF_i(X))^2,
\end{equation}
where $\omega^k=(\omega^k_1, \omega^k_2, \dots, \omega^k_n)$ is the weight vector at the $k$-th iteration, and $\odot$ denotes the element-wise product. Hence, we perform the line search at the $k$-th Newton iteration based on the equation \eqref{equation:3.1}. In this context, the line search criterion \eqref{equation:2.7} can be rewritten as follows:
\begin{equation}\label{equation:3.2}
f^k(X^k+\lambda^kS^k)\leq f^k(X^k)+\alpha \lambda^{k}\nabla f^k(X^k)^TS^k,
\end{equation}
where $f^{k}(\cdot)$ is given in \eqref{equation:3.1}.
If we analytically compute the expression for \(\nabla f^k(X^k)\) and substitute it into (\ref{equation:3.2}), we can further simplify the expression (\ref{equation:2.7}) to be
\begin{equation}\label{equation:3.3}
f^k(X^k+\lambda^kS^k)\leq f^k(X^k)+\alpha \lambda^{k} \left( \omega^k\odot(\omega^k\odot F(X^k)) \right)^T F^{\prime}(X^k)S^k.
\end{equation}

It is clear that $\alpha \nabla f^k(X^k)^TS^k$ remains in the form of a matrix-vector product, and the additional computational cost introduced by the equation (\ref{equation:3.3}) can be neglected. 
To balance the unbalanced nonlinearities in the nonlinear system using these weight multipliers, the update rule for the adaptive weight multiplier $\omega_{i}^k$ at the $k$-th iteration and the $i$-th component is given by the following formula:
\begin{equation}\label{equation:3.4}
\omega_{i}^{k+1} \leftarrow \delta^{k}_1\cdot \omega_{i}^{k} + \alpha^k  \left( \frac{|e_i^k|}{\Vert e^k\Vert_{\max}} + \delta^k_2\cdot\frac{\Vert e^k\Vert_{\max}-|e_i^k|}{\Vert e^k\Vert_{\max}} \right),  \quad i =1,2, \ldots, n,
\end{equation}
where $n$ denotes the total number of components in the nonlinear system of algebraic equations, $e^k:=F(X^k)$ denotes the residual vector at the $k$-th iteration, and $e^k_i:=F_i(X^k)$ represents the $i$-th component of the residual vector corresponding to the $k$-th iteration. 
Further, $\delta_1^k$ represents the weight decay factor, and $\delta_2^k$ represents the recognition factor, which are given by
\begin{align}\label{equation:3.5}
\delta_1^k=\delta \cdot \psi_1 \left( \frac{\Vert F(X^k) \Vert}{\Vert F(X^{k-1}) \Vert} \right), \quad
\delta_2^k=1-\psi_2 \left( \frac{\Vert F(X^k) \Vert}{\Vert F(X^{k-1}) \Vert} \right), 
\end{align}
where $\delta\in(0,1)$ denotes the base decay rate, and $\psi_1(t)=\exp ( -\frac{(t-1)^2}{2\sigma_1^2} )$ and $\psi_2(t)=\exp ( -\frac{(t-1)^2}{2\sigma_2^2} )$ are Gaussian functions with constants $\sigma_1$ and $\sigma_2$, respectively. 
Additionally, $\alpha^k$ in equation (\ref{equation:3.4}) is the learning rate for the $k$-th iteration, which is given by
\begin{equation}
\alpha^k=\alpha^{*}\cdot \left( \frac{2g^{k-1}}{g_{\max}} \right), \label{equation:3.6}
\end{equation}
where $\alpha^{*}$ denotes the initial learning rate, $g_{k-1}$ represents the number of line search iterations in the previous Newton iteration, and $g_{\max}$ denotes the maximum number of line search iterations allowed.

The update rule for the adaptive weight multiplier \eqref{equation:3.4} may seem intricate at first glance, but the underlying concept is straightforward. 
Below we provide some heuristic explanations.

\begin{itemize}
	\item
We can interpret the update rule \eqref{equation:3.4} for the adaptive weight multiplier as a special form of gradient ascent. In fact, the update rule \eqref{equation:3.4} that we propose is a gradient-free update scheme, which updates the weight vector $\omega^k$ in the direction of increasing residual norm. Thus, we can see that component in the nonlinear system that contribute a larger proportion to the residual norm will be assigned a larger weight multiplier, and vice versa. 
The term corresponding to equation \eqref{equation:2.10} in the INB method can be written as
\begin{equation}\label{equation:3.7}
	\alpha \nabla f^k(X^k)^TS^k=-\alpha \nabla f^k(X^k)^TF^{\prime}(X^k)^{-1}F(X^k)=-\alpha \Vert \omega^k\odot F(X^k)\Vert^2=-\alpha\sum_{i=1}^n \left( \omega^k_iF_i(X^k) \right)^2.
\end{equation} 
Combining the analysis from Section \ref{sec:New} and equation \eqref{equation:3.4}, we can see that the line search based on the merit function $f^k$ defined in \eqref{equation:3.1} can select a larger and more appropriate step length. This ensures that the components in the nonlinear system that contribute a larger proportion to the residual norm are sufficiently reduced in the current iteration, while also balancing the nonlinearities in the system.

\item
More specifically, we can observe that
$\left( \frac{|e_i^k|}{\| e^k \|_{\max}} + \delta_2^k \cdot \frac{\| e^k \|_{\max} - |e_i^k|}{\| e^k \|_{\max}} \right) \in [0, 1]$, and $\omega_i^1\neq 0, \forall i \in \{1, \dots, n\}$. Therefore, by performing a simple recursion and scaling based on equation \eqref{equation:3.4}, the bounds can be given by
\begin{equation}\label{equation:3.8}
\omega^k_i\in \left[ 0,\, \omega_i^1+\frac{2\alpha^*}{1-\delta} \right].
\end{equation}

\item
The weight decay factor $\delta_1^k$ determines the contribution of $\omega^k$ to $\omega^{k + 1}$. As the number of iterations increases, the contribution of $\omega^k$ to subsequent weight vectors gradually decreases. 
From equation \eqref{equation:3.5}, we can see that when the Newton-type solver stagnates, i.e. $\frac{\Vert F(X^k) \Vert}{\Vert F(X^{k-1}) \Vert}\approx 1$, this decay effect diminishes, thereby gradually increasing the focus on the unbalanced nonlinearities in the nonlinear system. In fact, the recognition factor $\delta_2^k$ and the learning rate $\alpha_k$ play a similar role to the weight decay factor. When the Newton iteration stagnates, $\delta_2^k \approx 0$, and the learning rate $\alpha_k$ typically takes a larger value. At this point, the update rule for the weight multipliers approaches the following rule:
\begin{equation}\label{equation:3.9}
\omega_{i}^{k+1} \leftarrow \delta\cdot \omega_{i}^{k} + 2\alpha^* \frac{|e_i^k|}{\Vert e^k\Vert_{\max}},  \quad i =1, 2, \ldots, n.
\end{equation}
\end{itemize}

The proposed adaptive strategy exhibits a high degree of generality and can be seamlessly incorporated into any Newton-type solver, including inexact Newton methods and their various refined variants, such as the Multilayer Nonlinear Elimination Preconditioned Inexact Newton (MNEPIN) \cite{luo2020multilayer} and the Preconditioned Inexact Newton with Learning Capability ($\text{PIN}^{\mathcal L}$) \cite{luo2023preconditioned}. In the following, we illustrate, taking the Inexact Newton method with backtracking (INB) as an example, how our adaptive strategy ought to be applied across different Newton-type solvers. 

The novel algorithm outlined above is summarized in Algorithm \ref{alg:ARDN}.

\begin{algorithm}
\caption{ARDN: Adaptive Residual-Driven Newton Solver}
\label{alg:ARDN}
\begin{algorithmic}[1]
\STATE \textbf{Input:} $X^0$, $\eta_0$, $\beta$, $\gamma_a$, $\gamma_r$, $\alpha$, $\rho$, $\omega^0$,$\delta$, $\alpha^*$, $\sigma_1$, $\sigma_2$
\STATE \textbf{Output:} $X^*$
\STATE $k\leftarrow 0$ \hfill  $\triangleright$ Nonlinear iteration index
\STATE $f(X):=\frac{\Vert F(X)\Vert^2}{2}$ \hfill$\triangleright$ Define the merit function
\WHILE{$\Vert F(X^k)\Vert> \max\{\gamma_r\Vert F(X^0)\Vert, \gamma_a\}$}
\STATE // Compute the forcing term $\eta^k$ using equation (\ref{equation:2.5}).
\STATE Select $\eta^k\in[0, 1)$\
\STATE $\Delta\omega^k\leftarrow \alpha^k \left( \frac{|e^k|}{\Vert e^k\Vert_{\max}} + \delta^k_2\cdot\frac{\Vert e^k\Vert_{\max}-|e^k|}{\Vert e^k\Vert_{\max}} \right)$  \hfill  $\triangleright$ Compute the weight increment.
\STATE $\omega^k \leftarrow \delta_1^k\omega^k + \Delta \omega^k$ \hfill $\triangleright$ Update the weight vector
\STATE // Precondtioner $\mathcal{P}^k$
\STATE $S^k \leftarrow \mathrm{Krylov}(F^{\prime}(X^k),-F(X^k), \mathcal{P}^k)$\hfill$\triangleright$ Solve the linear system (\ref{equation:2.3}) (at the level $\eta^k$)
\STATE $f^k(X):=\frac{\Vert \omega^k \odot F(X)\Vert^2}{2}$ \hfill$\triangleright$ Define the merit function at the \( k \)-th step.
\FOR{$i \in \{0, \dots, g_{max}\}$}
    \STATE // Obtain the step length through backtracking line search.
    \IF{$f^k(X^k + \lambda_i S^k) \leq f^k(X^k) + \alpha \lambda_i \nabla f^k(X^k)^T S^k$}
        \STATE $\lambda^k\leftarrow \lambda_i$
        \STATE break
        \ELSE
        \STATE $\lambda_i \leftarrow \rho \lambda_i $
    \ENDIF
\ENDFOR
\STATE $k\leftarrow k+1$ \hfill  $\triangleright$ Update nonlinear iteration index
\STATE $X^k \leftarrow X^k+\lambda^kS^k$   \hfill  $\triangleright$ Update the approximate solution
\ENDWHILE
\STATE $X^*\leftarrow X^k$ 
\end{algorithmic}
\end{algorithm}

Based on the preceding discussion, we now offer several remarks to further elucidate the ARDN algorithm.

\begin{remark}
The ratio of residual norms $\frac{\Vert F(X^k) \Vert}{\Vert F(X^{k-1}) \Vert}$ and the number of line search iterations in the last Newton step $g_{k-1}$ are considered as by-products obtained from the classical INB method, without incurring additional computational costs. This represents a key advantage of the ARDN algorithm.
\end{remark}

\begin{remark}
	The ARDN algorithm identifies and balances the slow components in the residual space using an adaptive multiplier that is driven by the residuals, which inherently does not necessitate additional analysis of the nonlinear system or reliance on physical information underlying the nonlinear system.
\end{remark}

\begin{remark}
It's worth noting that the residual-driven adaptive strategy does not involve preconditioning of the nonlinear system. Its essence is to enhance the line search, enabling it to choose a step length that is beneficial for all components to decrease sufficiently. Furthermore, the residual-driven adaptive strategy can also be combined with some preconditioning algorithms \cite{cai2002nonlinearly, dolean2016nonlinear, liu2015field} to produce even better results.
\end{remark}

\subsection{Alternative Adaptive Strategies for the ARDN algorithm}
In this section, we briefly discuss some alternative approaches for the ARDN algorithm. 
Specifically, the content to be discussed herein will serve as substitutes for the aforementioned update rule \eqref{equation:3.4}. Meanwhile, a brief exploration of some other potential avenues for improvement will also be presented.

\begin{itemize}
	\item \emph{Simplified ARDN Strategy I.} In the previous section, we introduced the ARDN algorithm and its main idea. Here, we present a simplified version of ARDN. In fact, we can replace the aforementioned update rule \eqref{equation:3.4} with the following expression:
	\begin{equation}\label{equation:3.10}
		\omega_{i}^{k+1} \leftarrow \delta \cdot\omega_{i}^{k} + \alpha^*  \frac{|e_i^k|}{\Vert e^k\Vert_{\max}} ,  \quad i =1, 2, \ldots, n.
	\end{equation}
	Update rule \eqref{equation:3.10} actually represents the original idea of our work, which is to continuously adjust the current merit function based on the direction of increasing residuals. Additionally, update rule \eqref{equation:3.10} can also be viewed as a special case of update rule \eqref{equation:3.4}. This is because when the standard deviations $\sigma_1\rightarrow\infty$ and $\sigma_2\rightarrow\infty$, the Gaussian function degrades into a constant function that always takes the value of 1. By setting $\alpha^k\equiv \alpha^*$, update rules \eqref{equation:3.4} and \eqref{equation:3.10} become equivalent in this scenario.
	
	\item \emph{Simplified ARDN Strategy II.} Building upon the aforementioned Simplified ARDN Strategy I, we can propose the following Simplified ARDN Strategy II:
	\begin{equation}\label{equation:3.11}
		\omega_{i}^{k+1} \leftarrow \delta \cdot\omega_{i}^{k} + \alpha^k  \frac{|e_i^k|}{\Vert e^k\Vert_{\max}} ,  \quad i =1, 2, \ldots, n,
	\end{equation}
	where $\alpha^k=\alpha^{*}\cdot \left( \frac{2g^{k-1}}{g_{\max}} \right)$. In fact, this is the only modification compared to update rule \eqref{equation:3.10} in update rule \eqref{equation:3.11}, where the fixed learning rate $\alpha^*$ is replaced by an adaptive learning rate $\alpha^k$. The purpose of this approach is to differentiate the importance of different iteration points and further improve the algorithm's performance. At this point, we can also multiply the previous weight accumulation term $\omega^{k-1}$ by the Gaussian function $\psi_1$ to weaken the influence of previous points on the current point.
	
	\item \emph{Employing different Newton-type solvers.} The ARDN algorithm presented in the previous section is an enhancement based on the Inexact Newton method with backtracking (INB). In fact, as we previously mentioned, we have also applied this strategy to other Newton-type solvers, such as  $\text{PIN}^{\mathcal L}$ \cite{luo2023preconditioned}, ASPIN \cite{cai2002nonlinearly}, and others.
	
	\item \emph{Improving Inexact Newton Directions.} Consider an optimization problem with a merit function $f = \frac{\|\omega^k\odot F\|^2}{2}$, where solving nonlinear system of algebraic equations can be transformed into solving the aforementioned optimization problem. For an adaptive merit function $f = \frac{\|\omega^k\odot F\|^2}{2}$, we can consider employing quasi-Newton methods or other techniques to solve the optimization problem, with the aim of improving the direction.
\end{itemize}

\section{Additional Optimizations in the Newton Solver}
\label{sec:addopt}
In this section, we recall a specific scheme intended to enhance the performance of the Newton-type solver, and we elaborate on how it is integrated with our adaptive strategy. This scheme will be applied in the numerical experiments presented in subsequent sections.

\subsection{Preconditioned inexact Newton methods with learning capability}
Preconditioned inexact Newton methods with learning capability are also referred to as $\text{PIN}^{\mathcal L}$ \cite{luo2023preconditioned}. 
This introduces an unsupervised learning strategy based on principal component analysis (PCA) \cite{chatterjee2000introduction}, designed to accelerate the performance of a Newton solver's algorithm. Its primary focus is on learning the detrimental behavior of the Newton solver in the residual subspace based on the training problem.
By solving a nonlinear preconditioning system based on principal component analysis, a new initial guess is obtained. In this section, we will briefly review the $\text{PIN}^{\mathcal L}$ algorithm and elucidate how it integrates with our adaptive strategy. 

Considering the nonlinear system of algebraic equations \eqref{equation:2.1}. 
In order to construct a preconditioner using PCA, it is necessary to collect the data set first. 
By iterating for $s-1$ ($s\geq 1$) steps with the INB algorithm, we obtain a set of approximate solution sequences $\{X^k\}_{k=0}^{s-1}$ and the corresponding residual sequences $\{ F(X^k) \}_{k=0}^{s-1}$. Define that
\[\bar{F}=\frac1s\sum_{k=0}^{s-1}F(X^k),\, \bar{X}=\frac1s\sum_{k=0}^{s-1}X^k \text{  and  }
\hat{F}^k=F(X^k)-\bar{F},\, \hat{X}^k=X^k-\bar{X}.
\]
We centralize the collected data set and form the centralized residual matrix $R$ and solution matrix $S$, where 
\begin{align}
R= \left[ \hat{F}^0, \hat{F}^1, \dots, \hat{F}^{s-1} \right] \in \mathbb{R}^{n \times s}, \quad 
S = \left[ \hat{X}^0, \hat{X}^1, \dots, \hat{X}^{s-1} \right] \in \mathbb{R}^{n \times s}. \label{equation:4.1}
\end{align}
Subsequently, we apply PCA to the residual matrix $R$ and the solution matrix $S$ respectively. 

In fact, we need to solve the following two optimization problems to obtain the projection operator $P\in \mathbb{R}^{n \times d}$ and  $Q\in \mathbb{R}^{n \times d}$,
\begin{align}
\max_{P\in H_{n\times d}}\mathcal{V}_1(P),\quad
\max_{Q\in H_{n\times d}}\mathcal{V}_2(Q), \label{equation:4.2}
\end{align}
where $H_{n\times d}=\begin{Bmatrix} A|A\in\mathbb{R}^{n\times d},A^TA=I_{d\times d} \end{Bmatrix}$, and
\begin{align} 
\mathcal{V}_1(P)=\sum_{k=0}^{s-1}\left\|P^T\hat{F}^k\right\|^2, \quad
\mathcal{V}_2(Q)=\sum_{k=0}^{s-1}\left\|Q^T\hat{X}^k\right\|^2,\label{equation:4.3}
\end{align}
where we specify the number of principal components as $d$. Generally speaking, the singular value decomposition (SVD) is used to solve the optimization problems given in \eqref{equation:4.2}, and thereby the projection operators $P$ and $Q$ can be obtained. 
We perform the SVD on the residual matrix $R$ and the solution matrix $S$ respectively as follows: 
\begin{align}
R=\hat{U}_R\hat{\Sigma}_R\hat{V}_R^T,\quad
S=\hat{U}_S\hat{\Sigma}_S\hat{V}_S^T, \label{equation:4.4}
\end{align}
where $\hat{U}_R$ and $\hat{U}_S$ are both $n\times n$ orthogonal matrices, $\hat{\Sigma}_R$ and $\hat{\Sigma}_S$ are $n\times s$ diagonal matrices with singular values in descending order, and $\hat{V}_R$ and $\hat{V}_S$ are $s\times s$ orthogonal matrices. The solutions of the optimization problems (\ref{equation:4.2}) are given by $P=\hat{U}_R^d$ and $Q=\hat{U}_S^d$, 
where $\hat{U}_R^d$ and $\hat{U}_S^d$ are composed of the first d columns of $\hat{U}_R$ and $\hat{U}_S$ respectively. 

Therefore, we can use the projection operator $P$ to construct the approximate nonlinear system. 
\begin{align}\label{equation:4.5}
    \mathcal{F}(Y):=PP^T \left( F(Y)-\bar{F} \right)+\bar{F}.
\end{align}
In order to solve the approximate nonlinear system \eqref{equation:4.5} by the Newton method, we need to solve the following linear system at each iteration.
\begin{align}\label{equation:4.6}
    PP^TF^{\prime}(Y^j)S^j=-\mathcal{F}(Y^j)=-PP^T \left( F(Y^j)-\bar{F} \right)-\bar{F},
\end{align}
where \(Y^j\) represents the current approximate solution, and \(S^j \in \mathbb{R}^n\) indicates the Newton direction in the j-th iteration. We multiply both sides of equation \eqref{equation:4.6} on the left by the matrix \(P^T\) and use the projection operator \(Q\), then we can obtain the equivalent form of the linear system \eqref{equation:4.6}. 
\begin{align}\label{equation:4.7}
    P^TF^{\prime}(Y^j) \left(QS^j_p \right)=-P^TF(Y^j),
\end{align}
where $S^j=QS^j_p$ and we denote \(\mathcal{F}_P(Y^j)=P^T\mathcal{F}(Y^j)=P^TF(Y^j)\) as the projected approximate nonlinear system, \(\mathcal{J}_p(Y^j)=P^TF^{\prime}(Y^j)Q\) as the projected Jacobian matrix, and \(S^j_p\) as the low-dimensional Newton correction. Therefore, we obtain the Jacobian system of the projected subspace 
\begin{align}\label{equation:4.8}
    \mathcal{J}_p(Y^j)S^j_p=-\mathcal{F}_p(Y^j).
\end{align}
We only need to solve the low-dimensional linear system \eqref{equation:4.8} instead of \eqref{equation:4.6}. The new approximate solution is given by 
\begin{align}\label{equation:4.9}
    Y^{j+1}=Y^j+QS^j_p.
\end{align}
And we adopt the following stopping criterion to determine whether the approximate nonlinear system converges,
\begin{align}\label{equation:4.10}
	\|\mathcal{F}(Y^*)\|\leq\gamma_r^s\|\mathcal{F}(Y^0)\|,
\end{align}
where $\gamma_r^s$ is the relative tolerance for the nonlinear system in the subspace. Subsequently, the approximate solution $Y^*$ obtained from solving the nonlinear system in the subspace can be used as the new initial guess for the INB algorithm. 
\begin{algorithm}
\caption{$\text{PIN}^{\mathcal L}$: Preconditioned inexact Newton methods with learning capability}
\label{alg:PINL}
\begin{algorithmic}[1]
\STATE \textbf{Input:} $X^0$, $\eta_0$, $\beta$, $\gamma_a$, $\gamma_r$, $\alpha$, $\rho$, $\gamma_r^s$
\STATE \textbf{Output:} $X^*$
\STATE // The INB algorithm yields the centered data matrices $R\in\mathbb{R}^{n \times s}$ and $S\in\mathbb{R}^{n \times s}$
\STATE $R, S = \mathrm{INB}(X^0, s)$  \hfill $\triangleright$ Obtain the centered data matrix
\STATE // Perform PCA on matrices $R$ and $S$ separately
\STATE $P,\, Q = \mathrm{PCA}(R, d),\,  \mathrm{PCA}(S, d)$ \hfill  $\triangleright$ Obtain the projection operators $P$ and $Q$.
\STATE $j\leftarrow 0$ \hfill  $\triangleright$ Subspace iteration index
\WHILE{$\|\mathcal{F}(Y^j)\|>\gamma_r^s\|\mathcal{F}(Y^0)\|$}
\STATE // Obtain a new initial guess $Y^*$ through subspace iteration
\STATE $\mathcal{F}_p\leftarrow P^TF(Y^j)$
\STATE $\mathcal{J}_p\leftarrow P^TF(Y^j)Q$ \hfill  $\triangleright$ Compute the Jacobian matrix
\STATE $\mathcal{J}_pS^j_p=\mathcal{F}_p$  \hfill  $\triangleright$ Solve the linear system accurately to obtain the subspace Newton direction
\STATE $j \leftarrow j+1$
\STATE $Y^j \leftarrow Y^j + QS^j_p$  \hfill  $\triangleright$ Update
\ENDWHILE
\STATE $Y^* \leftarrow Y^j$ \hfill  $\triangleright$ a new initial guess
\STATE Use the new initial guess as the initial value for the INB algorithm.
\STATE $X^* \leftarrow \mathrm{INB}(Y^*)$ \hfill  $\triangleright$ Compute the approximate solution using the INB algorithm
\end{algorithmic}
\end{algorithm}
We can summarize the above ideas into the following Algorithm \ref{alg:PINL}:

To explain how the $\text{PIN}^{\mathcal L}$ algorithm integrates with our method, we provide the following remarks:

\begin{remark}
	The $\text{PIN}^{\mathcal L}$ algorithm firstly employs the classic INB algorithm for iteration, generating a training set from which a projection operator is obtained. Subsequently, this projection operator is used to precondition the original system. By solving the nonlinear system in the subspace, a good initial guess is obtained, which serves as the good initial guess for the global Newton iteration.
\end{remark}

\begin{remark}
	Besides using the INB algorithm on the original problem to obtain the training set, we can also apply the INB algorithm to a training problem to obtain the training set. Here, the training problem refers to a simpler but similar problem compared to the original problem
\end{remark}

\begin{remark}
	Besides using the INB algorithm to obtain the training set, other Newton-type solvers can also be used to generate the training set. Therefore, for the $\text{PIN}^{\mathcal L}$ algorithm, integrating it with the ARDN algorithm can be achieved by simply replacing the classic INB method in both the training step and the global INB step with the ARDN algorithm
\end{remark}

\section{Numerical experiments}
\label{sec:exam}
In this section, we conducted numerical experiments on a series of nonlinear system of algebraic equations to evaluate the performance of the algorithm we proposed. 
In the numerical experiments, each linear system obtained from each Newton iteration is solved using the GMRES method \cite{saad2003iterative}. The codes for the nonlinear solver and the linear solver are all implemented in Python. Among them, the numerical calculations of GMRES and the Jacobian matrix are implemented by calling SciPy \cite{gommers2022scipy}. The absolute tolerance and relative tolerance of the Newton-type solver are set to $10^{-8}$ and $10^{-12}$ respectively. The line search decay rate $\rho = 0.5$, and the initial forcing term $\eta_0 = 0.25$. The maximum number of nonlinear iterations is set to $200$. 
If it is greater than $200$, we consider that the method fails to be convergent and therefore stop iterating.

Here we fix some notions.

$N_{ite}$:  The total number of Newton iterations performed; 

$T(s)$:  The total running time of the algorithm with Unit seconds;

$N_{sta}$:  The number of stagnation times during the entire algorithm operation.

Here stagnation means that the remaining amount of the current step is almost the same as the remaining amount of the previous step,
that is, $\frac{\Vert F(X^k) \Vert}{\Vert F(X^{k-1}) \Vert}\approx 1$. However, in order to accurately describe the stagnation phenomenon, we stipulate that the current Newton iteration is considered to be in a stagnant state when the following conditions are met: 
\begin{equation}\label{eq:5.1}
\left|  \Vert F(X^{k})\Vert-\Vert F(X^{k-1})\Vert  \right| \leq 10^{-6} \Vert F(X^{k})\Vert.
\end{equation}
Unless otherwise specified, we always use \eqref{eq:5.1} to determine whether the Newton-type solver is stagnant. 

\subsection{The chemical equilibrium system}
Consider the chemical equilibrium system \cite{grosan2008new}:
\begin{equation}  
\left\{  
\begin{array}{l}  \label{equation:5.2}
x_1 x_2 + x_1 - 3x_5=0, \\  
2x_1 x_2 + x_1 + x_2 x_3^2 + R_8 x_2 - R x_5 + 2R_{10} x_2^2 + R_7 x_2 x_3 + R_9 x_2 x_4=0, \\  
2x_2 x_3^2 - 8x_5 + R_6 x_3 + R_7 x_2 x_3=0, \\  
R_9 x_2 x_4 + 2x_4^2 - 4R x_5=0, \\  
x_1(x_2 + 1) + R_{10} x_2^2 + R_8 x_2 + R_5 x_3^2 - 1 + R_6 x_3 + R_7 x_2 x_3 + R_9 x_2 x_4=0.  
\end{array}  
\right.  
\end{equation}
where $R = 10, R_5 = 0.193, 
R_6 = \frac{0.002597}{\sqrt{40}}, 
R_7 = \frac{0.003448}{\sqrt{40}}, 
R_8 = \frac{0.00001799}{40}, 
R_9 = \frac{0.0002155}{\sqrt{40}}, 
R_{10} = \frac{0.00003846}{40}$.
For this problem, we take the initial value $X^0 = (0,\dots, 0)\in \mathbb{R}^5$, and the Jacobian matrix is calculated through numerical approximation. 
Although this is a low-dimensional problem, the solution to this problem has certain difficulties. (i) The Jacobian matrix generated in each Newton iteration is extremely ill-conditioned. (ii) The INB algorithm applied to this problem is slowly convergent, and even the failure of convergence (we can observe this phenomenon through Figure \ref{fig:res} or Table \ref{tab:1} and Table \ref{tab:2} later).

\subsubsection{Algorithm Performance Evaluation}
We first verify the effectiveness of the ARDN Algorithm \ref{alg:ARDN} we proposed and compare its numerical results with the INB Algorithm \ref{alg:INM} and the $\text{PIN}^{\mathcal L}$ Algorithm \ref{alg:PINL}. Here, we set the maximum number of line search iterations $g_{max} = 36$, and the other parameters remain the same as before. 
In the $\text{PIN}^{\mathcal L}$ method, the number of training sets is $s = 8$, and the number of principal components is $d = 2$. 

\begin{figure}[h]
	\centering
	\includegraphics[width=0.8\textwidth]{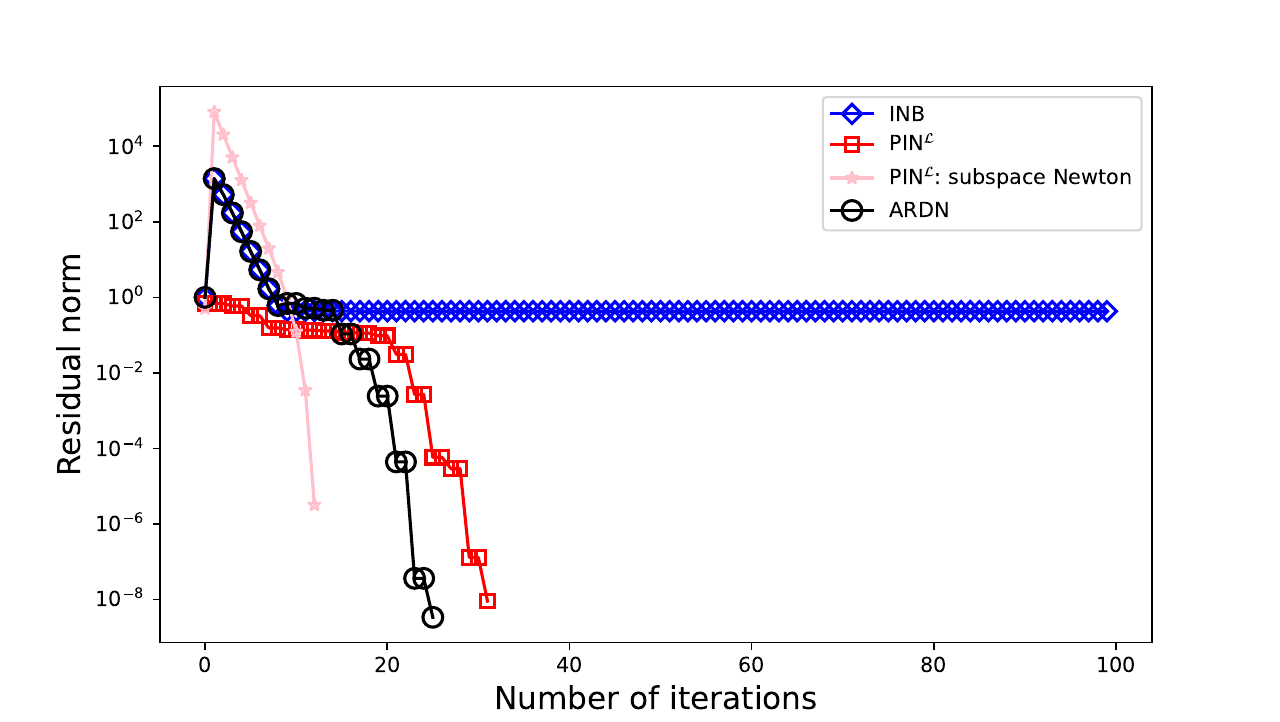}
	\caption{The nonlinear residual history obtained using the INB, ARDN and $\text{PIN}^{\mathcal L}$ methods}
	\label{fig:res}
\end{figure}

Figure \ref{fig:res} shows the evolution history of the nonlinear residuals obtained using the INB Algorithm \ref{alg:INM}, the ARDN Algorithm \ref{alg:ARDN}, and the $\text{PIN}^{\mathcal L}$ Algorithm \ref{alg:PINL}. Table \ref{tab:1} presents the relevant results on the performance indicators of the three algorithms. We can observe that the INB algorithm fails to converge even after reaching the maximum number of iterations of $N_{ite}=200$, the $\text{PIN}^{\mathcal L}$ algorithm converges after $N_{ite}=31$ iterations, and our proposed algorithm only needs $N_{ite}=25$ iterations to meet the convergence criteria. 
This not only confirms the effectiveness of our ARDN algorithm on this typical problem, but also has a faster convergence rate compared to existing INB algorithm and ARDN algorithm.

Additionally, we can observe that the three algorithms exhibit different degrees of stagnation. Among them, the INB algorithm stagnates for $N_{sta}=164$ steps, and it can be observed that after the nonlinear residual norm reaches $\mathcal{O}(10^{-1})$, it is unable to decrease and remains in a stagnant state, ultimately leading to the failure of convergence. The $\text{PIN}^{\mathcal L}$ algorithm has a slight stagnation phenomenon, with a total of $N_{sta}=15$ stagnant steps. Compared to the INB algorithm, it has a significant improvement. The ARDN algorithm only stagnates for $N_{sta}=8$ steps, generally showing a rapid convergence phenomenon. This result is in line with our expectations. This is because the algorithm we designed itself has certain corrections for the stagnation phenomenon. 

\begin{table}
\vspace*{0.2cm} \caption{The results obtained using IN, $\text{PIN}^{\mathcal L}$, and ARDN for the chemical equilibrium system.} \label{tab:1}
\begin{center}
\begin{tabular}{cccccccccc}
\hline & INB & \multicolumn{3}{c}{$\text{PIN}^{\mathcal L}$} & ARDN\\
\cline{3-5} 
& &  Training & Subspace Newton & Global IN & & \\
\hline
$N_{ite}$ & 200 & 8 & 12 & 31 & 25\\
$T(s)$ &0.1091 & 0.0035 & 0.0016 & 0.0172 & 0.0120\\
$N_{sta}$ &164 & --- & --- & 15 & 8\\
\hline
\end{tabular}
\end{center}
\end{table}

\subsubsection{Impact of preselected parameters on the algorithm}
To study the impact of parameters on the performance of the ARDN algorithm, we first examine the chemical equilibrium system using different values of $\alpha_*$ and $g_{max}$. The relevant numerical results are presented in Table \ref{tab:2}. We can observe that the algorithm has good robustness to the parameter $\alpha_*$, while the selection of $g_{max}$ has a significant impact on the algorithm. Generally speaking, the impact of the maximum number of iterations of the line search, $g_{max}$, on the algorithm is difficult to avoid. This is because the value of $g_{max}$ affects the size of the step length. If the value of $g_{max}$ is too small, it may cause the merit function to not decrease sufficiently; conversely, it may lead to a too small step length, which easily causes the algorithm to fall into a stagnant state. 

\begin{table}[htbp]
\vspace*{0.2cm} \caption{The impact of parameters $\alpha_*$ and $g_{max}$ on the performance of ARDN for the Chemical Equilibrium Application} \label{tab:2}
\begin{center}
\begin{tabular}{cccccccccccccccc}
\hline 
&\multicolumn{3}{c}{\textbf{$g_{max}=12$}} && \multicolumn{3}{c}{\textbf{$g_{max}=24$}} && \multicolumn{3}{c}{\textbf{$g_{max}=36$}} && \multicolumn{3}{c}{\textbf{$g_{max}=48$}} \\
\cline{2-4} \cline{6-8} \cline{10-12}  \cline{14-16}
\textbf{$\alpha_*$} & \textbf{$N_{ite}$} & \textbf{$T(s)$} & \textbf{$N_{sta}$} && 
    \textbf{$N_{ite}$} & \textbf{$T(s)$} & \textbf{$N_{sta}$} && \textbf{$N_{ite}$} & 
   \textbf{$T(s)$} & \textbf{$N_{sta}$} && 
   \textbf{$N_{ite}$} & \textbf{$T(s)$} &
   \textbf{$N_{sta}$}\\
\hline
   0.03 & 71 & 0.0278 & 14 && 55 & 0.0236 & 18 && 25 & 0.0125 & 8 && 36 & 0.0157 & 15 \\
   0.06 & 71 & 0.0288 & 14 && 55 & 0.0240 & 18 && 25 & 0.0126 & 8 && 36 & 0.0158 & 15 \\
   0.24 & 71 & 0.0287 & 14 && 55 & 0.0240 & 18 && 25 & 0.0125 & 8 && 36 & 0.0156 & 15 \\
   1.00 & 71 & 0.0288 & 14 && 55 & 0.0242 & 18 && 25 & 0.0126 & 8 && 36 & 0.0157 & 15 \\
\hline
\end{tabular}
\end{center}
\end{table}

The standard deviation $\sigma_1$ and $\sigma_2$ of the Gaussian function are another crucial factor that influences the algorithm's performance. For various standard deviations employed in ARDN, we present the comparative results in Table \ref{tab:3}. The magnitude of the standard deviation determines the \emph{width} of the Gaussian function, so the specific values of the standard deviations $\sigma_1$ and $\sigma_2$ are crucial. On the one hand, when the standard deviation $\sigma_2$ is very small, the ARDN algorithm becomes almost indistinguishable from the classical INB method, or at best shows only minor improvements. When $\sigma_2$ is extremely large, the ARDN algorithm resembles the alternative adaptive strategies for the ARDN algorithm mentioned previously. On the other hand, if the standard deviation $\sigma_1$ is small, the adaptive weight multipliers from previous Newton iterations generated by the ARDN algorithm have negligible influence on the current Newton iteration. Conversely, if $\sigma_1$ is large, the attenuation effect provided by the corresponding Gaussian function $\psi_1$ can be considered negligible, which may diminish the influence of the adaptive weight multipliers in the current Newton iteration step.
In both scenarios, the effectiveness of the ARDN algorithm is compromised. Therefore, selecting the standard deviations $\sigma_1$ and $\sigma_2$ within a suitable range is crucial. Generally, we recommend choosing both standard deviations, $\sigma_1$ and $\sigma_2$, to be appropriately small for two reasons: firstly, to mitigate the influence of previous weights on the current Newton iteration step, and secondly, to enable the algorithm to more accurately identify stagnation phenomena in the Newton solver.

\begin{table}[htbp]
\vspace*{0.2cm} \caption{The impact of parameters $\sigma_1$, $\sigma_2$ and $g_{max}$ on the performance of ARDN for the Chemical Equilibrium Application} \label{tab:3}
\begin{center}
\begin{tabular}{ccccccccccccc}
\hline 
\multicolumn{6}{c}{ARDN} && \multicolumn{3}{c}{INB}\\
\cline{1-6} \cline{8-10} 
\textbf{$g_{max}$} & \textbf{$\sigma_1$} & \textbf{$\sigma_2$} & \textbf{$N_{ite}$} & \textbf{$T(s)$} &\textbf{$N_{sta}$} && \textbf{$N_{ite}$} & \textbf{$T(s)$} & \textbf{$N_{sta}$}
\\
\hline
   12 & 0.30  & 0.01 & 99 & 0.0449 & 28 && 99 & 0.0449 &28 \\
   12 & 0.30  & 0.25 & 71 & 0.0284 & 14 && 99 & 0.0447 &28 \\
   12 & 0.30  & 0.50  & 71 & 0.0289 & 14 && 99 & 0.0449 &28 \\
   12 & 0.30  & 100  & 71 & 0.0285 & 14 && 99 & 0.0455 &28 \\
   12 & 0.01 & 0.25 & 71 & 0.0283 & 14 && 99 & 0.0445 &28 \\
   12 & 0.50  & 0.25 & 71 & 0.0289 & 14 && 99 & 0.0449 &28 \\
   12 & 100  & 0.25 & 99 & 0.0449 & 28 && 99 & 0.0455 &28 \\
\hline
   24 & 0.30  & 0.01 & 45 & 0.0216 & 13 && 61 & 0.0274 &21 \\
   24 & 0.30  & 0.25 & 45 & 0.0222 & 13 && 61 & 0.0275 &21 \\
   24 & 0.30  & 0.50  & 55 & 0.0237 & 18 && 61 & 0.0276 &21 \\
   24 & 0.30  & 100  & 69 & 0.0330 & 24 && 61 & 0.0271 &21 \\
   24 & 0.01 & 0.25 & 55 & 0.0227 & 18 && 61 & 0.0273 &21 \\
   24 & 0.50  & 0.25 & 45 & 0.0219 & 13 && 61 & 0.0280 &21 \\
   24 & 100  & 0.25 & 61 & 0.0278 & 21 && 61 & 0.0275  &21 \\
\hline
\end{tabular}
\end{center}
\end{table}

\subsubsection{Integration of ARDN with other acceleration strategies}
In this subsection, we once again emphasize that our adaptive weighting strategy can be applied to any Newton-type solver, including but not limited to the classic INB, $\text{PIN}^{\mathcal L}$ \cite{luo2023preconditioned}, ASPIN \cite{hwang2005parallel}, and other modified Newton solvers \cite{luo2020multilayer} \cite{novello2024accelerating}. Here, we only test how our adaptive enhancement strategy integrates with $\text{PIN}^{\mathcal L}$ and observe its performance. Considering the chemical equilibrium system, when the maximum number of line search $g_{max}$ is greater than $42$, the number of inexact Newton methods no longer changes. Therefore, we may as well set $g_{max} = 42$. At the same time, we reasonably choose the parameters of the $\text{PIN}^{\mathcal L}$ algorithm to optimize its performance. Specifically, we select the number of training sets $s=8$ and the number of principal components $d=2$. 

\begin{figure}[h]
	\centering
	\includegraphics[width=0.8\textwidth]{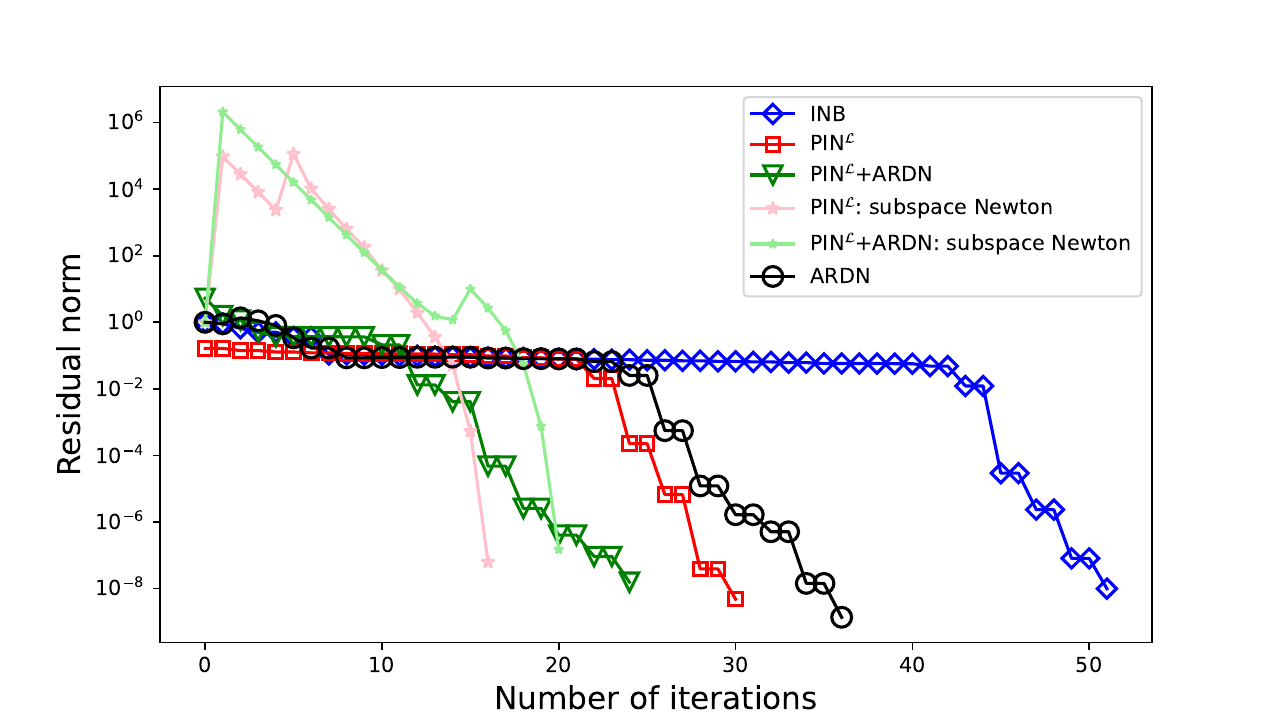}
	\caption{The nonlinear residual history obtained using the INB, ARDN, $\text{PIN}^{\mathcal L}$, and $\text{PIN}^{\mathcal L}$ + ARDN methods}
	\label{fig:compar}
\end{figure}

Figure \ref{fig:compar} exhibits the history of nonlinear residuals obtained using INB, $\text{PIN}^{\mathcal L}$, ARDN, and ARDN+$\text{PIN}^{\mathcal L}$ for a given scenario. We can find that the $\text{PIN}^{\mathcal L}$+ ARDN algorithm only iterates $24$ steps, while the $\text{PIN}^{\mathcal L}$ algorithm iterates $30$ steps, and the running time of the algorithm is also lower. Obviously, the algorithm performance has been improved. In fact, we can adjust the number of training sets $s$ and the number of principal components $d$, and the specific numerical results are presented in Table \ref{tab:4}. We can observe that for most cases, the algorithm performance of $\text{PIN}^{\mathcal L}$+ARDN is better than that of $\text{PIN}^{\mathcal L}$ algorithm, and the robustness of $\text{PIN}^{\mathcal L}$ algorithm is enhanced. Therefore, we have reasons to believe that the combination of $\text{PIN}^{\mathcal L}$ algorithm and ARDN algorithm can achieve better results. 

\begin{table}[htbp]
\vspace*{0.2cm} \caption{Numerical Results of $\text{PIN}^{\mathcal L}$ and $\text{PIN}^{\mathcal L}$ + ARDN methods for the Chemical Equilibrium System under different $d$ and $s$} \label{tab:4}
\begin{center}
\begin{tabular}{ccccccccccccccc}
\hline 
\multicolumn{5}{c}{$\text{PIN}^{\mathcal L}$} && \multicolumn{5}{c}{$\text{PIN}^{\mathcal L}$+ARDN}\\
\cline{1-5} \cline{7-11} 
\textbf{$s$} & \textbf{$d$} & \textbf{$N_{ite}$} & \textbf{$T(s)$} &\textbf{$N_{sta}$} && \textbf{$s$} & \textbf{$d$} & \textbf{$N_{ite}$} & \textbf{$T(s)$} &\textbf{$N_{sta}$}
\\
\hline
   6  & 2 & 200 & 0.140 & 191 && 6 & 2 & 40 & 0.023 & 18 \\
   7  & 1 & 151  & 0.075 &  74  && 7 & 1 & 45 & 0.028 & 21 \\
   7  & 2 & 200 & 0.140 & 192 && 7 & 2 & 23 & 0.020 & 11 \\
   8  & 1 & 175  & 0.082 &  86  && 8 & 1 & 36 & 0.022 & 16 \\
   8  & 2 & 30  & 0.030 &  15  && 8 & 2 & 24 & 0.017 & 10 \\
   9  & 1 & 151  & 0.077 &  74  && 9 & 1 & 22 & 0.016 & 9 \\
   10 & 1 & 149  & 0.067 &  73  && 10& 1 & 18 & 0.013 & 7 \\
   16 & 1 & 59   & 0.034 &  28  && 16& 1 & 17 & 0.015 & 6 \\
   16 & 2 & 30   & 0.024 &  15  && 16& 2 & 19 & 0.021 & 9 \\
\hline
\end{tabular}
\end{center}
\end{table}
\subsection{Convection-Diffusion Problem}
Consider the following nonlinear partial differential equation
\begin{equation}  
\left\{  
\begin{array}{l}  \label{equation:5.4}
\Delta u(x,y)+Cu(x,y)(u_x(x,y)+u_y(x,y))=g(x,y), (x,y)\in\Omega,\\
u(x,y)=0, (x,y)\in\partial\Omega,
\end{array}  
\right.  
\end{equation}
where $C$ is a constant and $\Omega=[0,1]\times[0,1]$. This problem is given by \cite[page 109]{kelley1995iterative}. The function $g$ is selected such that the solution of the equation is the following function 
\begin{align}\label{eq:44}
    u(x,y)=10xy(1-x)(1-y)\exp(x^{4.5}).
\end{align}

For this problem, we use the standard central difference scheme for discretization. After discretization, a nonlinear system of algebraic equations can be obtained. When $C$ is relatively small, it is relatively simple to solve the nonlinear system of algebraic equations obtained by discretizing the nonlinear partial differential equation \eqref{equation:5.4}; however, as $C$ gradually increases, the difficulty of solving the nonlinear system of equations will continue to increase. We will continuously increase the value of $C$ and adjust the grid size to observe the performance of different methods in solving the problem. The specific numerical results are presented in Table \ref{tab:5}. For this problem, we take $X^0 = (0,\dots,0)$ as the initial value for iteration. we replace $10^{-6}$ in equation \eqref{eq:5.1} with $10^{-2}$, set $g_{max} = 24$, $\sigma_1 = 0.3$, $\sigma_2 = 0.25$, and set the GMRES restart number to $50$. And we set the absolute tolerance and the relative tolerance to $10^{-8}$ and $10^{-10}$, respectively. Moreover, we do not adopt additional preconditioning techniques, and the other parameters are the same as before. 
\begin{table}[htbp]
\vspace*{0.2cm} \caption{Performance Indicators for Convection-Diffusion Problems with Different Parameters $C$ and Grid Sizes} \label{tab:5}
\begin{center}
\begin{tabular}{ccccccccccc}
\hline 
&& \multicolumn{3}{c}{INB} && \multicolumn{3}{c}{ARDN}\\
\cline{3-5} \cline{7-9} 
\textbf{$C$} & \textbf{$n\times n$} & \textbf{$N_{ite}$} & \textbf{$T(s)$} &\textbf{$N_{sta}$} && \textbf{$N_{ite}$} & \textbf{$T(s)$} &\textbf{$N_{sta}$}
\\
\hline
   80  & $50\times 50$   & 60  & 0.558  & 10  && 35 & 0.286 & 2 \\
   80  & $100\times 100$ & 67  & 3.662  & 13  && 47 & 2.753 & 5 \\
   100 & $50\times 50$   & 121 & 1.309  & 74  && 73 & 0.742 & 4 \\
   100 & $100\times 100$ & 171 & 14.08 & 140 && 87 & 7.197 & 27 \\
   120 & $50\times 50$   &200 &2.790  &164 &&108 & 1.341 & 6 \\
   120 & $100\times 100$ & 191 &599.0 & 139 && 78 & 5.319 & 5 \\
   140 & $50\times 50$   &200 &3.124  & 176&& 138 & 1.810 & 17 \\
   140 & $100\times 100$ &200 &17.52 & 187&& 84 & 6.114 & 7 \\
\hline
\end{tabular}
\end{center}
\end{table}

From Table \ref{tab:5}, we can observe that with the increase of $C$ or the number of grids, the number of iterations required by the INB method generally shows an upward trend, and when $C > 120$, the INB method basically fails to converge. Additionally, we can observe that for all possible values of $C$ and grid sizes, the number of iterations required by the proposed ARDN method is significantly smaller than that required by the INB method, and the running time of the ARDN method is also more advantageous. Moreover, the number of stagnation of the ARDN method is significantly reduced compared to the INB method, with almost no more than $20$. Thus, we can find that the ARDN method effectively alleviates the stagnation phenomenon in the INB method and accelerates the convergence of the solver.

For this problem, we also studied the integration of the ARDN algorithm and the $\text{PIN}^{\mathcal L}$ algorithm. We uniformly set the number of training sets $s = 6$ and the number of principal components $d = 3$. The numerical results are presented in Table \ref{tab:6}. From Table \ref{tab:6}, we can observe that, in most cases, the effect of the $\text{PIN}^{\mathcal L}$ + ARDN algorithm is better than that of the $\text{PIN}^{\mathcal L}$ algorithm. Therefore, we have reason to believe that the integration of the ARDN algorithm and the $\text{PIN}^{\mathcal L}$ algorithm can achieve better results. 

\begin{table}[htbp]
\vspace*{0.2cm} \caption{Numerical Results of $\text{PIN}^{\mathcal L}$ and $\text{PIN}^{\mathcal L}$ + ARDN Methods under Different C and Grid Sizes} \label{tab:6}
\begin{center}
\begin{tabular}{ccccccccccccccc}
\hline 
\multicolumn{5}{c}{$\text{PIN}^{\mathcal L}$} && \multicolumn{5}{c}{$\text{PIN}^{\mathcal L}$+ARDN}\\
\cline{1-5} \cline{7-11} 
\textbf{$C$} & \textbf{$n\times n$} & \textbf{$N_{ite}$} & \textbf{$T(s)$} &\textbf{$N_{sta}$} && \textbf{$C$} & \textbf{$n\times n$} & \textbf{$N_{ite}$} & \textbf{$T(s)$} &\textbf{$N_{sta}$}
\\
\hline
80  & $50\times 50$  & 71 & 0.763 & 18 && 80 & $50\times 50$ & 41 & 0.420 & 4 \\
80  & $100\times 100$  & 69 & 5.037 & 18 && 80 & $100\times 100$ & 77 & 5.926 & 28 \\
90  & $50\times 50$  & 91 & 1.109 & 36 && 90 & $50\times 50$ & 60 & 0.774 & 6 \\
100  & $50\times 50$ & 104  & 1.218 &  70  && 100 & $50\times 50$ & 73 & 0.841 & 15 \\
100  & $100\times 100$ & 113  & 7.647 & 80  && 100 & $100\times 100$ & 69 & 5.635 & 15 \\
110  & $50\times 50$ & 134  & 1.861 &  102  && 110 & $50\times 50$ & 102 & 1.340 & 20 \\
110  & $100\times 100$ & 150  & 12.56 &  126  && 110 & $100\times 100$ & 53 & 4.140 & 11 \\
140  & $50\times 50$ & 200  & 3.148 &  190  && 140 & $50\times 50$ & 143 & 2.038 & 37 \\
\hline
\end{tabular}
\end{center}
\end{table}

\subsection{A Series of Challenging Problems}
In this subsection, we select a series of challenging problems Problem 1 (P1) to Problem 5 (P5). For the majority of these problems, either the classical Newton method fails to converge or the INB method converges relatively slowly. We conducted tests on our proposed algorithm concerning these problems, while also comparing it with the INB method. In the following problem description, mod means the remainder after integer division.

\subsubsection{Problem 1}
Modified Rosenbrock \cite{friedlander1997solving}.
Let $n$ be an even and non-zero natural number, and 
\begin{align*}
F_{k}(x) & =\frac{1}{1+\exp(-x_{k})}-0.73, &&\quad\mathrm{mod}(k,2)=1, \notag\\
F_{k}(x) & =10(x_{k}-x_{k-1}^2), &&\quad\mathrm{mod}(k,2)=0,\quad 
k=1,\dots, n. \notag
\end{align*}
We use $X^0=(X_1,\ldots,X_n)^T$ as the initial value, where $X_{l}=-1.8$ for $\mathrm{mod}(l,2)=1$ and $X_{l}=-1$ for  $\mathrm{mod}(l,2)=0$, 
This problem shows a slow convergence when using the INB method (this phenomenon can be observed from Table \ref{tab:7}). 

\subsubsection{Problem 2} Augmented Rosenbrock \cite{friedlander1997solving}.
Let $n$ be a multiple of 4 and a non-zero natural number and
\begin{align*}
& F_{k}(x)=10(x_{k+1}-x_{k}^{2})&&\quad\mathrm{mod}(k,4)=1, \notag\\
& F_{k}(x)=1-x_{k-1}&&\quad\mathrm{mod}(k,4)=2, \notag\\
& F_{k}(x)=1.25x_{k}-0.25x_{k}^{3}&& \quad\mathrm{mod}(k,4)=3,\notag\\
& F_{k}(x)=x_{k} &&\quad\mathrm{mod}(k,4)=0, \quad k=1,\dots,n. \notag
\end{align*}
We use $X^0=(X_1,\ldots,X_n)^T$ as the initial value, where $X_{l}=-1.2$ for $\mathrm{mod}(l,4)=1$, $X_{l}=1$ for $\mathrm{mod}(l,4)=2$, $X_{l}=-1$ for $\mathrm{mod}(l,4)=3$ and $X_{l}=20$ for $\mathrm{mod}(l,4)=4$.
This problem will lead to the convergence failure when using the classical Newton method \cite{friedlander1997solving}.

\subsubsection{Problem 3} Tridiagonal \cite{lukvsan1994inexact}. Let $n$ be a non-zero natural number and 
\begin{align*}
& F_k(x)=4(x_k-x_{k+1}^2), & & k=1, \notag\\
& F_k(x)=8x_k(x_k^2-x_{k-1})-2(1-x_k)+4(x_k-x_{k+1}^2), & & k=2,\ldots,n-1, \notag\\
& F_k(x)=8x_k(x_k^2-x_{k-1})-2(1-x_k), & & k=n. \notag
\end{align*}
We use $X^0=(12,\ldots,12)^T\in\mathbb{R}^n$ as the initial value.
This problem converges slowly when the INB method is applied (this phenomenon can be observed from Table \ref{tab:7}).

\subsubsection{Problem 4} Five-diaggonal \cite{lukvsan1994inexact}. 
Let $n$ be a non-zero natural number and
\begin{align*}
& F_k(x)=4(x_1-x_2^2)+x_2-x_3^2,& & k=1,  \notag\\
& F_k(x)=8x_2(x_2^2-x_1)-2(1-x_2)+4(x_2-x_3^2)+x_3-x_4^2, & & k=2, \notag\\
& F_k(x)=8x_k(x_k^2-x_{k-1})-2(1-x_k)+4(x_k-x_{k+1}^2)+x_{k-1}^2-x_{k-2}+x_{k-1}-x_{k-2}^2,&& k=3,\ldots,n-2, \notag\\
& F_{k}(x)=8x_{n-1}(x_{n-1}^2-x_{n-2})-2(1-x_{n-1})+4(x_{n-1}-x_n^2)+x_{n-2}^2-x_{n-3},& & k=n-1,  \notag\\
& F_k(x)=8x_n(x_n^2-x_{n-1})-2(1-x_n)+x_{n-1}^2-x_{n-2},& & k=n. \notag 
\end{align*}
We use $X^0=(12,\ldots,12)^T\in\mathbb{R}^n$ as the initial value. 
This problem converges slowly when the INB method is applied (this phenomenon can be observed from Table \ref{tab:7}).

\subsubsection{Problem 5} Tridimensional valley \cite{friedlander1997solving}.
Let $n$ be a multiple of 4 and a non-zero natural number and
\begin{align*}
F_{k}(x) & =(c_{2}x_{k}^{3}+c_{1}x_{k})\exp\left(\frac{-x_{k}^{2}}{100}\right)-1&&\quad\mathrm{mod}(k,3)=1, \notag\\
F_{k}(x) & =10(\sin(x_{k-1})-x_{k}),&&\quad\mathrm{mod}(k,3)=2, \notag\\
F_{k}(x) & =10(\cos(x_{k-2})-x_{k}), &&\quad\mathrm{mod}(k,3)=0, \quad k=1,\ldots,n, \notag
\end{align*}
where $c_1=1.003344481605351, c_2=-3.344481605351171\times10^{-3}$.  We use $X^0=(X_1,\ldots,X_n)^T\in\mathbb{R}^n$ as the initial value,
 where $X_{l}=-4$ for $\mathrm{mod}(l,3)=1$,  $X_{l}=1$ for $\mathrm{mod}(l,3)=2$, $X_{l}=2$ for $\mathrm{mod}(l,3)=0$.
This problem is not convergent for the classical Newton method \cite{friedlander1997solving}.

We applied the INB algorithm and the ARDN algorithm to the problems P1-P5. For each test problem, we selected different dimensions $n$ for experiments. We conducted multiple tests and counted the number of stagnations $N_{sta}$ of the INB algorithm and the ARDN algorithm respectively, as well as the total number of Newton iterations $N_{ite}$ of the corresponding problems, the algorithm running time $T(s)$, and other performance indicators. 
These numerical results are presented in Table \ref{tab:7}. The Jacobian matrices of all problems are analytically calculated. The GMRES restart number is set to $50$, 
$g_{max} = 12$. For the linear equation systems in each Newton iteration, we did not adopt preconditioning. For P1 and  P5, we change $10^{-6}$ in equation \eqref{eq:5.1} to $10^{-2}$; For problems P2, P3 and P4, we make no changes. Here, we did not adjust the values of the standard deviations $\sigma_1$ and $\sigma_2$, and they were fixed at $\sigma_1 = 0.3$ and $\sigma_2 = 0.25$ for all problems, and the other parameters are the same as before. 

\begin{table}[htbp]
\vspace*{0.2cm} \caption{The numerical results of problems P1 - P5 obtained by the ARDN method and the INB method} \label{tab:7}
\begin{center}
\begin{tabular}{ccccccccccc}
\hline 
& \multicolumn{3}{c}{INB} && \multicolumn{3}{c}{ARDN}\\
\cline{2-4} \cline{6-8} 
 \textbf{Problems} & \textbf{$N_{ite}$} & \textbf{$T(s)$} &\textbf{$N_{sta}$} && \textbf{$N_{ite}$} & \textbf{$T(s)$} &\textbf{$N_{sta}$}
\\
\hline
P1,\, $n=6.0 \times 10^{1}$  & 80 & 0.0296 & 27 && 55 & 0.0155  & 17 \\
P1,\, $n=6.0 \times 10^{2}$  & 66 & 0.0233 & 15 && 54 & 0.0200  & 12  \\
P1,\, $n=6.0 \times 10^{3}$  & 65 & 0.0669 & 12  && 55 & 0.0555  & 9  \\
\hline 
P2,\, $n=6.0 \times 10^{3}$  & 11 & 0.0065 & 1  && 10 & 0.0045 & 0  \\
P2,\, $n=4.0 \times 10^{5}$  & 13 & 1.0308 & 1  && 10 & 0.7054  & 0 \\
P2,\, $n=8.0 \times 10^{7}$  & 13 & 26.514 & 1  && 10 & 17.293  & 0 \\
\hline 
P3,\, $n=6.0 \times 10^{1}$  & 77 & 3.3662  & 28 && 53 & 1.6998  & 14  \\
P3,\, $n=1.2 \times 10^{3}$  & 90 & 131.53  & 34 && 60 & 44.679  & 11 \\
P3,\, $n=2.4 \times 10^{3}$  & 96 & 522.13  & 39 && 60 & 187.42  & 10 \\
\hline 
P4,\, $n=1.0 \times 10^{2}$  & 117 & 10.032 & 51  && 49 & 1.5824 & 8 \\
P4,\, $n=1.0 \times 10^{3}$  & 98 & 177.55 & 39  && 62 & 56.007 & 12 \\
P4,\, $n=4.0 \times 10^{3}$  & 98 & 1465.5 & 37 && 58 & 276.89 & 7  \\
\hline 
P5,\, $n=1.2 \times 10^{3}$  & 39 & 0.1778 & 5  && 32 & 0.1411  & 3 \\
P5,\, $n=4.8 \times 10^{3}$  & 39 & 0.4037 & 5  && 32 & 0.2433  & 3  \\
P5,\, $n=9.6 \times 10^{3}$  & 39 & 0.8523 & 5  && 32 & 0.7220  & 3  \\
\hline
\end{tabular}
\end{center}
\end{table}

From Table \ref{tab:7}, it is evident that the ARDN method consistently outperforms the INB method across all tested problems. This superiority is reflected in both the total number of Newton iterations $N_{ite}$ and the algorithm's running time $T(s)$, where the ARDN method demonstrates notable improvements over the INB method.
Moreover, a closer examination of the number of stagnations $N_{sta}$ reveals a significant trend. For every problem, the ARDN method exhibits a reduction in the number of stagnations compared to the INB method, albeit to varying degrees. This observation aligns with the findings from our earlier experiments. The decrease in stagnation events plays a crucial role in enhancing the convergence rate of the Newton method. Specifically, the greater the reduction in stagnations, the more pronounced the decrease in the number of nonlinear iterations required by the Newton method, thereby contributing to overall computational efficiency.

\section{Concluding remarks}
\label{sec:con}
We have proposed an adaptive residual-driven Newton-type solver for solving nonlinear system of algebraic equations. It assigns a weight multiplier to each component in the nonlinear system of algebraic equations, and these weight multipliers are adaptively updated based on the residuals. The main purpose is to identify and balance the strong nonlinearity in the nonlinear system. The main features of the method are that it can dynamically identify and balance the nonlinearity of the system, and it has the characteristics of low cost and effectively alleviating the stagnation phenomenon in the Newton-type solver. 
Numerical results show that our algorithm is more robust and converges more quickly compared to the INB algorithm and the $\text{PIN}^{\mathcal L}$ algorithm. Furthermore, we point out that the ARDN algorithm can be easily combined with other preconditioning techniques or other methods, and has the potential to yield better robustness and faster convergence. Strict convergence rate analysis is an interesting topic, and it will be reported in a future work.

\bibliographystyle{alpha}
\bibliography{ARDN}

\end{document}